\documentclass{article}

\usepackage{makeidx}
\usepackage{multirow}
\usepackage{epsf}
\usepackage{latexsym}
\usepackage{natbib}
\usepackage{tabularx}
\usepackage{multicol}
\usepackage{comment}

\usepackage{amsmath,amsfonts,amscd,amsxtra,amssymb}
\usepackage{array,eqnarray}
\usepackage{graphicx,epsfig,lscape,color}
\usepackage{natbib}
\bibpunct{(}{)}{;}{a}{,}{,}

\usepackage{multirow}
\usepackage{lineno}

\DeclareMathOperator{\E}{E}				

\begin{document}
\title{A note on Tempelmeier's $\beta$-service measure under non-stationary stochastic demand}
\author{Roberto Rossi, Onur A. Kilic, S. Armagan Tarim}

\maketitle

\begin{abstract}
\citet{Tempelmeier2007} considers the problem of computing replenishment cycle policy parameters under non-stationary stochastic demand and service level constraints. He analyses two possible service level measures: the minimum no stock-out probability per period ($\alpha$-service level) and the so called ``fill rate'', that is the fraction of demand satisfied immediately from stock on hand ($\beta$-service level). For each of these possible measures, he presents a mixed integer programming (MIP) model to determine the optimal replenishment cycles and corresponding order-up-to levels minimizing the expected total setup and holding costs. His approach is essentially based on imposing service level dependent lower bounds on cycle order-up-to levels. In this note, we argue that Tempelmeier's strategy, in the $\beta$-service level case, while being an interesting option for practitioners, does not comply with the standard definition of ``fill rate''. By means of a simple numerical example we demonstrate that, as a consequence, his formulation might yield sub-optimal policies.
\end{abstract}

\section{Introduction}

The increasing pace in new product developments has resulted in shorter product life-cycles through which demands do not follow stationary patterns \citep{Kurawarwala1996,Graves2008}. Henceforth, inventory problems addressing non-stationary stochastic demands have gained a growing interest from both researchers and practitioners \citep{citeulike:7928534}. A recent paper, \citet{Tempelmeier2007}, addresses the replenishment cycle policies under non-stationary stochastic demand and $\beta$-service level (i.e., ``fill rate'') constraints. The \textit{fill rate} is the fraction of demand satisfied immediately from stock on hand. $\beta$-service level constraints therefore specify the minimum prescribed fraction of customer demand that should be met routinely, without backorders or lost sales. 

Tempelmeier's work constitutes an interesting development that extends results such as those presented by \citet{citeulike:7766622} to the non-stationary stochastic demand case. More specifically, Tempelmeier extends the model proposed in \citet{Tarim2004}, by replacing the  $\alpha$-service level constraints --- which enforce a minimum no-stockout probability per period --- with a new set of constraints based on the inverse first-order loss function. It is stated that the resultant formulation provides the optimal replenishment cycle plan under $\beta$-service level constraints. In this research note, we argue that the $\beta$-service level formulation proposed in \citet{Tempelmeier2007} does not comply with the standard definition of $\beta$-service level found in the literature.

In what follows, we provide the formal definition of $\beta$-service level (Section \ref{sec:definition}) and we discuss its application in an inventory system controlled with a replenishment cycle policy. Then we discuss the formulation proposed in \citeauthor{Tempelmeier2007} (Section \ref{sec:tempelmeier}) for computing optimal replenishment cycle policy parameters under $\beta$-service level constraints and show, by means of a simple numerical example, that this formulation may yield suboptimal policy parameters (Section \ref{sec:example}).

\section{$\beta$-service measure}\label{sec:definition}

The $\beta$-service level is a well established service measure used in many practical applications and has been covered by many textbooks on inventory control \citep[see e.g.][]{Silver1998,Axsater2006}. \cite{Axsater2006} defines $\beta$-service level as the fraction of demand satisfied immediately from stock on hand. This definition is formalized within the context of finite horizon inventory models as follows \citep[see e.g.][]{Chen2003,Thomas2005}:
\begin{equation}\label{beta}
1 - \E\left\{\frac{\text{Total backorders within the planning horizon}}{\text{Total demand within the planning horizon}}\right\}.
\end{equation}

The replenishment cycle policy divides the finite planning horizon into a number of, say $m$, consecutive replenishment cycles. We can re-write (\ref{beta}) by taking those into account as
\begin{equation}\label{beta_cycle}
1 - \E\left\{\frac{\sum_{i=1}^m\text{Total backorders within the $i$'th replenishment cycle}}{\sum_{i=1}^m\text{Total demand within the $i$'th replenishment cycle}}\right\}.
\end{equation}

\section{\citeauthor{Tempelmeier2007}'s formulation} \label{sec:tempelmeier}

For the ease of exposition, here we only provide the formulation of the $\beta$-service level constraints. The reader is referred to \cite{Tempelmeier2007} for the rest of the model. The set of constraints proposed by \citeauthor{Tempelmeier2007} to impose $\beta$-service level are as follows:

\begin{equation}\label{cons_tempelmeier}
\E\{I_t\}\geq\sum_{j=1}^t\left[F^{-1}_{Y^{(t-j+1,t)}}(\beta)-\sum_{i=t-j+1}^{t}\E\{D_i\}\right]P_{tj},\quad t=1,\ldots,T
\end{equation}

\noindent where, $I_t$ is the net inventory position at the end of period $t$; $F^{-1}_{Y^{(t-j+1,t)}}$ is the inverse loss function of the total demand in periods $(t-j+1,\ldots,t)$; $D_t$ is the random demand in period $t$; and, $P_{tj}$ is the binary indicator variable that is equal to 1 if the last replenishment before period $t$ takes place in period $t-j+1$ and to 0 otherwise. Following Tempelmeier, the expected net inventory position is assumed to be non-negative; however, it should be noted that relaxing non-negativity constraints on expected net inventory positions may yield better $\beta$-service plans in terms of expected cost. This issue is beyond the scope of this note and therefore not addressed here.

Eq.(\ref{cons_tempelmeier}) is binding only if the indicator variable $P_{tj}$ is equal to 1. Let us consider a replenishment cycle covering periods $(t'-j'+1,\ldots,t')$, i.e. $P_{t'j'}=1$. Then the binding part of the constraint reads:

\begin{equation}\label{cons_tempelmeier2}
\E\{I_{t'}\} + \sum_{i=t'-j'+1}^{t'}\E\{D_i\} \geq F^{-1}_{Y^{(t'-j'+1,t')}}(\beta).
\end{equation}

The replenishment at period $t'-j'+1$ covers the interval $(t'-j'+1,\ldots,t')$. The left hand side of the inequality represents the order-up-to level for period $t'-j'+1$. The constraint clearly imposes a lower bound on the order-up-to level for this cycle. Therefore, when these constraints are used, the same $\beta$-service level is imposed on each and every cycle within the planning horizon. This corresponds to the following definition of $\beta$-service level:

\begin{equation}\label{beta_tempelmeier}
1 - \max_{i=1,\ldots,m}\left[\E\left\{\frac{\text{Total backorders in replenishment cycle $i$}}{\text{Total demand in replenishment cycle $i$}}\right\}\right].
\end{equation}

It is clear that Eq.(\ref{beta_cycle}) is different from Eq.(\ref{beta_tempelmeier}). The original definition imposes a $\beta$-service level throughout the whole planning horizon, whereas \citeauthor{Tempelmeier2007}'s definition imposes a $\beta$-service level on each replenishment cycle within the planning horizon independently. The main difference is that, the former allows the decision maker to have $\beta$-service levels smaller than the specified level for individual cycles, while guaranteeing the specified level for the whole of the planning horizon, whereas the latter guarantees the specified $\beta$-service level for each replenishment cycle. It should be noted that Tempelmeier's strategy may be favorable for practitioners, since it allows a better control of the fill-rate provided to customers in each cycle. In practice, enforcing a given fill rate over
the whole planning horizon, rather than on each cycle separately, guarantees a lower cost at the expense of a varying individual replenishment cycle fill rates. Managers may therefore be interested in paying an additional price in order to have a better control over the fill rate provided in each cycle. For a thorough discussion on theoretical vs versus applied models in inventory control \citep[see][]{citeulike:8061205}.

\vspace{-1em}
\section{A numerical example}\label{sec:example}

Let us now consider a limit situation in which we aim to compute optimal non-stationary $(R,S)$ policy parameters for a 2-period planning horizon. The fixed ordering cost is 0, implying that the optimal plan has a replenishment in each period. The holding cost is 1. Period demands are normally distributed $N(\mu,\sigma)$ with parameters $N_1(1000,200)$ and $N_2(2000,200)$. We enforce a $\beta$-service level constraint with $\beta=0.98$.

According to \cite{Tempelmeier2007}, the minimum expected buffer stock level for period 1 is $181$ units (corresponding to an order-up-to-level, $R_1$, of 1181 units), which guarantees a $\beta$-service level of exactly $0.98$. Furthermore the minimum expected buffer stock level for period 2 is 99 units (corresponding to an order-up-to-level, $R_2$, of 2099 units), which guarantees a $\beta$-service level of exactly $0.98$.

From the first-order loss function, it is easy to see that the expected amount of items backordered in periods 1 and 2 are $19.90$ and $39.86$, respectively. It follows that the overall fill rate is $[(1000+2000)-(19.90+39.86)]/(1000+2000)=0.98$. The expected total holding cost of this solution is $280$ (i.e., 181 for period 1 and 99 for period 2).

A fill rate of 0.98 can be achieved also with the following plan, which guarantees a lower expected total holding cost.

For the first replenishment cycle, we fix an expected buffer stock level of 171, giving an order-up-to-level of $1171$. This expected buffer stock level is lower than the minimum expected buffer stock level allowed in Tempelmeier's model. In fact, it only guarantees, for the first replenishment cycle, a fill rate of 0.96 and an expected number of backorders of 21.79.

For the second replenishment cycle, we target a buffer stock level of 105, giving an order-up-to-level of 2105. This expected buffer stock level is higher than the minimum expected buffer stock level allowed in Tempelmeier's model. It guarantees, for the second replenishment cycle, a fill rate equal to 0.98 and an expected number of backordered items equal to 38.03.

In this alternative plan the overall $\beta$-service level over the 2-period planning horizon is $[(1000+2000)-(21.79+38.03)]/(1000+2000)=0.98$ as targeted. However, the required service level is attained with a lower expected total holding cost of 276 (i.e., 171 for period 1 and 105 for period 2).

%

\bibliographystyle{abbrvnat}
\bibliography{note}

\end{document}